\documentclass[12pt]{article}

\usepackage{amsfonts}
\usepackage{mathtools}
\usepackage{amsthm}
\usepackage{graphicx}
\usepackage[margin=1.5in]{geometry}
\usepackage{fancyhdr}
\usepackage[notquote]{hanging}
\usepackage{setspace}
\usepackage{MnSymbol}
\usepackage[hang,flushmargin]{footmisc}

\makeatletter
\providecommand*{\cupdot}{%
  \mathbin{%
    \mathpalette\@cupdot{}%
  }%
}
\newcommand*{\@cupdot}[2]{%
  \ooalign{%
    $\m@th#1\cup$\cr
    \sbox0{$#1\cup$}%
    \dimen@=\ht0 %
    \sbox0{$\m@th#1\cdot$}%
    \advance\dimen@ by -\ht0 %
    \dimen@=.5\dimen@
    \hidewidth\raise\dimen@\box0\hidewidth
  }%
}
\makeatother

\makeatletter
\def\old@comma{,}
\catcode`\,=13
\def,{%
  \ifmmode%
    \old@comma\discretionary{}{}{}%
  \else%
    \old@comma%
  \fi%
}
\makeatother

\title{Delta-matroids whose twist polynomials are monomials}
\author{Daniel Yuschak}
\date{}

\begin{document}
\sloppy

\maketitle

\begin{center}
\textit{Department of Mathematics, The Ohio State University, 231 West 18th Avenue, Columbus, OH 43210}
\end{center}

\noindent
\rule{\linewidth}{.5pt}

\kern6pt
\noindent
\textbf{Abstract}

\kern6pt
\noindent
The twist polynomial of a delta-matroid was recently introduced by Yan and Jin, who proved a characterization of binary delta-matroids whose twist polynomials are monomials. In this paper, we extend this result to all delta-matroids by proving that any delta-matroid whose twist polynomial is a monomial must be binary.\\



\noindent
\rule{\linewidth}{.5pt}\\

\let\thefootnote\relax\footnotetext{Email address: danielyuschak@gmail.com\\Declarations of interest: none}

\kern4pt
\noindent
\textbf{1. Introduction}\\

The concept of partial duality was introduced by Chmutov in [7] as a generalization of the notion of the dual of a ribbon graph. Gross, Mansour and Tucker [12] used this concept to define the partial duality polynomial, encoding the Euler genera of all partial duals of a ribbon graph. They also conjectured that there is no nontrivial orientable ribbon graph whose partial duality polynomial is a monomial. An infinite family of counterexamples to this conjecture was found by Yan and Jin in [15], and these have been shown to be essentially the only counterexamples [8, 16]. Chmutov and Vignes-Tourneret [8] then raised the possibility of generalizing the scope of this problem from ribbon graphs to delta-matroids, which are a generalization of matroids where bases may be of varying cardinality and which are closely connected to the theory of ribbon graphs. The delta-matroid analogue of the partial duality polynomial, the twist polynomial (see Definition 11 below), was defined by Yan and Jin in [17], and they then characterized all binary delta-matroids whose twist polynomials are monomials:\\

\noindent
\textbf{Theorem 1.} [18] \textit{The twist polynomial of a normal binary delta-matroid is a monomial if and only if each connected component of the looped simple graph corresponding to the delta-matroid is either a complete graph of odd order or a single vertex with a loop.}\\

Binary delta-matroids constitute a vanishing proportion of all delta-matroids. Indeed, consider all delta-matroids on a ground set of size $n$, equivalenced up to isomorphism and twisting. The number of nonequivalent binary delta-matroids is at most $2^{n(n+1)/2}$, the number of symmetric $n \times n$ binary matrices, while it follows from Theorem 1.1 of [10] that the total number of nonequivalent delta-matroids grows doubly exponentially in $n$. This trend can be seen for small $n$ in the following table:

\begin{center}
\begin{tabular}{ |c|c|c| }
\hline
Ground & Number of nonequivalent & Number of nonequivalent \\
set size & binary delta-matroids & delta-matroids \\
\hline
2 & 5 & 5 \\
\hline
3 & 13 & 16 \\
\hline
4 & 40 & 90 \\
\hline
5 & 141 & 2902 \\
\hline
\end{tabular}
\end{center}

\noindent
(The code used for this table and a list of the delta-matroids found are available as ancillary files at https://arxiv.org/abs/2208.13258.)

In this paper, we complete the classification of delta-matroids whose twist polynomials are monomials by showing that none of the twist polynomials of these nonbinary delta-matroids are monomials:\\

\noindent
\textbf{Theorem 2.} \textit{If the twist polynomial of a delta-matroid is a monomial, then the delta-matroid is binary.}\\

Using Theorems 1 and 2, we describe all delta-matroids whose twist polynomials are monomials in Section 3. Theorem 2 is proven by first characterizing the normal minors of delta-matroids as the restrictions of their normal twists, which allows several intermediate claims within the proof of Theorem 1 to be strengthened into necessary excluded-minor conditions for the class of delta-matroids whose twist polynomials are monomials. This then allows the use of the excluded-minor characterization of binary delta-matroids proven by Bouchet and Duchamp in [5] to show that the delta-matroid must be binary.\\

\noindent
\textbf{2. Definitions and Basic Properties}\\

We first review several concepts in delta-matroid theory relevant to the aim of this paper. Some of the terms defined here have analogous definitions in matroid theory, in which case the definitions here recover the usual ones when restricted to matroids.

A \textit{set system} is a pair $(E, \mathcal{F})$ consisting of a finite set $E$ called the \textit{ground set} along with a collection $\mathcal{F}$ of subsets of $E$ called \textit{feasible sets}.\\

\noindent
\textbf{Definition 3.} [2] A \textit{delta-matroid} is a set system $(E, \mathcal{F})$ such that $\mathcal{F}$ is nonempty and satisfies the \textit{symmetric exchange axiom}: For any $X, Y \in \mathcal{F}$ and any $u \in X \Delta Y$, there exists a $v \in X \Delta Y$ such that $X \Delta \{u, v\} \in \mathcal{F}$. Note that $v = u$ is permitted.\\

A delta-matroid is called \textit{normal} when the empty set is feasible, and a delta-matroid with all sets in $\mathcal{F}$ of the same size is a \textit{matroid} (described by its bases). We will denote the power set of a set $E$ by $\mathbb{P}(E)$. \\

\noindent
\textbf{Definition 4.} [14] The \textit{rank function} $r_M: \mathbb{P}(E) \mapsto \mathbb{N}_0$ of a matroid $M = (E, \mathcal{F})$ is the function given by
\[r_M(X) \coloneq \max_{F \in \mathcal{F}}|X \cap F| \textrm{ for all } X \in \mathbb{P}(E).\]

Additionally, for a matroid $M$ the expression $r(M) \coloneq r_M(E)$ will denote the size of any feasible set (i.e., basis) of $M$, which is the \textit{rank} of $M$.

For a delta-matroid $D = (E, \mathcal{F})$, $\mathcal{F}_{\min+i}$ will denote the collection of sets in $\mathcal{F}$ of cardinality $\min_{F \in \mathcal{F}} |F| + i$, and $\mathcal{F}_{\max-i}$ will denote those of cardinality $\max_{F \in \mathcal{F}} |F| - i$. In general, these collections do not necessarily satisfy the symmetric exchange axiom, however this is the case when $i = 0$ [4, Section 4]. Since all sets in these collections are the same size, these form matroids when $i = 0$, which will be denoted by $D_{\min} = (E, \mathcal{F}_{\min})$ and $D_{\max} = (E, \mathcal{F}_{\max})$.

There are various operations that can be performed on set systems, most importantly twists, deletions, and contractions. Let $D = (E, \mathcal{F})$ be a set system, and let $e$ be an element of $E$.\\

\noindent
\textbf{Definition 5.} [2] The \textit{twist} of $D$ by $e$ is the set system $D*e \coloneq (E, \{F \Delta \{e\}: F \in \mathcal{F}\})$.\\

\noindent
\textbf{Definition 6.} [9] The set system obtained from $D$ by \textit{deleting} $e$ is $D \backslash e \coloneq (E \backslash \{e\}, \{F: e \notin F, F \in \mathcal{F}\})$ if there exists $F \in \mathcal{F}$ with $e \notin F$, and $D \backslash e \coloneq (E \backslash \{e\}, \{F \backslash \{e\}: F \in \mathcal{F}\})$ otherwise.\\

\noindent
\textbf{Definition 7.} [9] The set system obtained from $D$ by \textit{contracting} $e$ is $D/e \coloneq (E \backslash \{e\}, \{F \backslash \{e\}: e \in F, F \in \mathcal{F}\})$ if there exists $F \in \mathcal{F}$ with $e \in F$, and $D/e \coloneq (E \backslash \{e\}, \mathcal{F})$ otherwise.\\

If $D$ is a delta-matroid, then all of these set systems are also delta-matroids [2, Remark 2] [5, Section 2].

It is straightforward to check that for any set system, twists commute with twists, deletions, and contractions when the elements of $E$ involved are distinct, and $D*e \backslash e = D/e$ and $D*e/e = D \backslash e$ [9, Lemma 2.7]. For set systems, deletions and contractions generally do not commute with other deletions and contractions, but it can be shown that they do commute when $D$ is a delta-matroid [6, Remark 20]. As such, when $D$ is a delta-matroid and $A$ is any subset of $E$, the delta-matroids $D*A$, $D \backslash A$, and $D/A$ are defined respectively by twisting, deleting, and contracting by the elements of $A$ in any order. The delta-matroid $D \backslash A^c$ obtained by deleting all elements except those in $A$ will be notated as $D|_A$ and is called the \textit{restriction} of $D$ to $A$. The delta-matroid $D*E$ is called the \textit{dual} of $D$ and denoted by $D^*$.\\

\noindent
\textbf{Definition 8.} A delta-matroid $D^\prime$ is a \textit{minor} of a delta-matroid $D$ if $D^\prime$ can be obtained from $D$ via some sequence of twists, deletions, and contractions.\\

We will make use of one more operation on delta-matroids, namely, the direct sum of two delta-matroids.\\

\noindent
\textbf{Definition 9.} [11] Given delta-matroids $D_1 = (E_1, \mathcal{F}_1)$ and $D_2 = (E_2, \mathcal{F}_2)$, the \textit{direct sum} of $D_1$ and $D_2$ is the delta-matroid $D_1 \oplus D_2 \coloneq (E_1 \cupdot E_2, \{F \subseteq E_1 \cupdot E_2: F \cap E_1 \in \mathcal{F}_1, F \cap E_2 \in \mathcal{F}_2\})$.\\

Delta-matroids naturally appear in ribbon graph theory, though only certain delta-matroids arise in this way. A slightly broader class of delta-matroids than those arising from ribbon graphs are binary delta-matroids, defined as follows.\\

\noindent
\textbf{Definition 10.} Let $E$ be a finite set and let $A = (a_{e, f})_{e, f \in E}$ be a symmetric $E \times E$ matrix over $\mathbb{F}_2$. For each $X \subseteq E$, let $A[X]$ denote the principal submatrix $(a_{e, f})_{e, f \in X}$. Then $(E, \{X \subseteq E: A[X]$ is nonsingular$\})$ (where $A[\O]$ is taken to be nonsingular) is a delta-matroid [3]. A \textit{binary} delta-matroid is any delta-matroid $D$ such that some twist of $D$ corresponds in this way to a symmetric binary matrix.\\

Symmetric $E \times E$ binary matrices are adjacency matrices of looped simple graphs with vertices bijectively labelled by $E$. Namely, a vertex labelled $e$ is looped if and only if $a_{e, e} = 1$, and vertices labelled $e$ and $f$ are connected by an edge if and only if $a_{e, f} = a_{f, e} = 1$. In this way, binary delta-matroids can equivalently be thought of as arising from looped simple graphs.

Given symmetric $E_i \times E_i$ binary matrices $A_i$ for $i \in \{1, 2\}$, each principal submatrix $(A_1 \oplus A_2)[X]$ of $A_1 \oplus A_2$ is equal to $A_1[X \cap E_1] \oplus A_2[X \cap E_2]$, and so $(A_1 \oplus A_2)[X]$ is nonsingular if and only if $A_1[X \cap E_1]$ and $A_2[X \cap E_2]$ are nonsingular. It thus follows that the direct sum of the delta-matroids arising from $A_1$ and $A_2$ is binary with the corresponding matrix $A_1 \oplus A_2$. In terms of looped simple graphs, this means that direct sums of binary delta-matroids correspond to disjoint unions of their underlying graphs.

The \textit{width} of a delta-matroid $D$ is the difference in size between its largest and smallest feasible sets, $r(D_{\max}) - r(D_{\min})$, and will be denoted by $w(D)$.\\

\noindent
\textbf{Definition 11.} [17] The \textit{twist polynomial} $^{\partial}w_D$ of a delta-matroid $D = (E, \mathcal{F})$ encodes the widths of its various twists, and is defined as follows:
\[^{\partial}w_D(z) \coloneq \sum_{A \subseteq E} z^{w(D*A)}\]

The twist polynomial of the direct sum of two delta-matroids is the product of their twist polynomials [17, Proposition 7], and so $^{\partial}w_{D_1 \oplus D_2}$ is a monomial if both $^{\partial}w_{D_1}$ and $^{\partial}w_{D_2}$ are.
\\\\
\noindent
\textbf{3. Examples (Counterexamples to the Gross-Mansour-Tucker Conjecture for Delta-Matroids)}\\

Theorems 1 and 2 together provide a complete characterization of delta-matroids whose twist polynomials are monomials in terms of looped simple graphs. In this section, we explicitly describe these delta-matroids in terms of their feasible sets.

Examples of delta-matroids whose twist polynomials are monomials are given below:\\

$D_3 \coloneq (\{1, 2, 3\}, \{\O, \{1, 2\}, \{1, 3\}, \{2, 3\}\})$

$D_5 \coloneq (\{1, 2, 3, 4, 5\}, \{\O, \{1, 2\}, \{1, 3\}, \{1, 4\}, \{1, 5\}, \{2, 3\}, \{2, 4\}, \{2, 5\}, \{3, 4\}, \{3, 5\}, \{4, 5\}, \{1, 2, 3, 4\}, \{1, 2, 3, 5\}, \{1, 2, 4, 5\}, \{1, 3, 4, 5\}, \{2, 3, 4, 5\}\})$\\

\noindent
In general, $D_{2k+1}$ is a delta-matroid on a ground set of cardinality $2k+1$ that is represented by a complete graph on $2k+1$ vertices, and its feasible sets are precisely the even-sized subsets of the ground set. Twisting $D_{2k+1}$ by any set of even size has no effect on the collection of feasible sets, and twisting it by any set of odd size results in the feasible sets being precisely the odd-sized subsets of the ground set. In either case, the width of the twisted delta-matroid is $2k$, and so the twist polynomial is given by $^{\partial}w_{D_{2k+1}}(z) = 2^{2k+1}z^{2k}$.

Since the twist polynomial is multiplicative under direct sums of delta-matroids, the twist polynomials of all finite direct sums $D = \bigoplus_i D_{2k_i+1}$ are monomials. These are the delta-matroids represented by disjoint unions of complete graphs of odd order, and constitute all counterexamples to the Gross-Mansour-Tucker conjecture among even normal delta-matroids.

Beyond even delta-matroids, there is another series of counterexamples $N_n$, where the ground set is of size $n$ and all subsets of the ground set are feasible. In fact, $N_n$ is the direct sum of $n$ copies of the delta-matroid $N_1 \coloneq (\{1\}, \{\O, \{1\}\})$, and so is represented by the disjoint union of $n$ looped vertices. It has width $w(N_n) = n$, and it is entirely unchanged by twisting. Thus its twist polynomial is the monomial $^{\partial}w_{N_n}(z) = 2^n z^n$.

Finally, it follows from Theorems 1 and 2 that the normal delta-matroids whose twist polynomials are monomials are precisely those of the form $D \oplus N_n$ for some even normal delta-matroid $D$ with $^{\partial}w_D$ a monomial and some $n \geq 0$. In terms of the feasible sets of delta-matroids, this means that for a normal delta-matroid $D = (E, \mathcal{F})$, $^{\partial}w_D$ is a monomial precisely when there is a partition of the ground set $E = A \cupdot B_1 \cupdot ... \cupdot B_m$ (where $A$ may be empty and $m$ may be 0) such that each $B_i$ is of odd size, and such that a set $F \subseteq E$ is feasible if and only if the numbers $|F \cap B_i|$ are all even. Note that for any delta-matroid of this form, twisting by any feasible set doesn't alter the delta-matroid.\\

\noindent
\textbf{4. Proof of the Main Result}\\

\noindent
\textbf{Lemma 12.} \textit{Suppose $D = (E, \mathcal{F})$ is a delta-matroid with a minor $M$. Then there exist $F, A \subseteq E$ such that $M = (D*F) \backslash A$ where $A$ is disjoint from some feasible set of $D*F$. If, moreover, $M$ is normal, then $F \in \mathcal{F}$.}\\

\noindent
\textbf{Proof.} Choose a sequence of twists, deletions, and contractions that transforms $D$ into $M$. Recall that for any delta-matroid $D'$ and any elements $e, f$ of the ground set, $D'/e = D'*e \backslash e$, and if $D' \backslash e*f$ is defined, $D' \backslash e*f = D'*f \backslash e$. Thus this sequence can be transformed into one only involving twists and deletions, and it can then be rearranged to be of the form $M = (D*X) \backslash A$ for some $X, A \subseteq E$. Choose a partition $A = Y \cupdot Z$ of $A$ such that $Y$ is a maximal subset of $A$ disjoint from at least one feasible set of $D*X$. Then the feasible sets of $(D*X) \backslash Y$ are the feasible sets of $D*X$ that are disjoint from $Y$, each of which contains $Z$. Therefore $M = (D*X) \backslash A = ((D*X) \backslash Y) \backslash Z = ((D*X) \backslash Y*Z) \backslash Z = ((D*X*Z) \backslash Y) \backslash Z = (D*F) \backslash A$ where $F = X \Delta Z$. Now for any feasible set $F'$ of $(D*X) \backslash Y$, $F' \Delta Z = F' \backslash Z$ is a feasible set of $D*F$ that is disjoint from $A$.

Since $A$ is disjoint from some feasible set of $D*F$, the feasible sets of $M$ form a subcollection of the collection of feasible sets of $D*F$. As such, if $M$ is normal, so is $D*F$, in which case $F \in \mathcal{F}$. \qedsymbol\\

\noindent
\textbf{Proposition 13.} [1, Lemma 1.2] \textit{For any delta-matroid $D = (E, \mathcal{F})$ and any $F_0 \in \mathcal{F}$, there exist $F_1 \in \mathcal{F}_{\min}$ and $F_2 \in \mathcal{F}_{\max}$ such that $F_1 \subseteq F_0 \subseteq F_2$.}\\

\noindent
\textbf{Theorem 14.} [5] \textit{A delta-matroid is binary if and only if it has no minor isomorphic to one of the following delta-matroids:}

\textit{1. $(\{1, 2, 3\}, \{\O, \{1, 2\}, \{1, 3\}, \{2, 3\}, \{1, 2, 3\}\})$}

\textit{2. $(\{1, 2, 3\}, \{\O, \{1\}, \{2\}, \{3\}, \{1, 2\}, \{1, 3\}, \{2, 3\}\})$}

\textit{3. $(\{1, 2, 3\}, \{\O, \{2\}, \{3\}, \{1, 2\}, \{1, 3\}, \{1, 2, 3\}\})$}

\textit{4. $(\{1, 2, 3, 4\}, \{\O, \{1, 2\}, \{1, 3\}, \{1, 4\}, \{2, 3\}, \{2, 4\}, \{3, 4\}\})$}

\textit{5. $(\{1, 2, 3, 4\}, \{\O, \{1, 2\}, \{1, 4\}, \{2, 3\}, \{3, 4\}, \{1, 2, 3, 4\}\})$}\\

\noindent
For the following results, we fix two small delta-matroids:\\

$D_1 \coloneq (\{1, 2\}, \{\O, \{1\}, \{2\}\})$

$D_2 \coloneq (\{1, 2, 3\}, \{\O, \{1, 2\}, \{1, 3\}\})$\\

\noindent
Note that this is not the same $D_1$ as in Section 3.\\

\noindent
\textbf{Lemma 15.} Any nonbinary delta-matroid has a minor isomorphic to $D_1$ or $D_2$.\\

\noindent
\textbf{Proof.} By Theorem 14, it suffices to show that each of the five delta-matroids listed in that theorem has a minor isomorphic to $D_1$ or $D_2$.

For the first of these, $(\{1, 2, 3\}, \{\O, \{1, 2\}, \{1, 3\}, \{2, 3\}, \{1, 2, 3\}\})*\{1, 2, 3\}|_{\{1, 2\}} = (\{1, 2\}, \{\O, \{1\}, \{2\}\}) \cong D_1$.

For the second, $(\{1, 2, 3\}, \{\O, \{1\}, \{2\}, \{3\}, \{1, 2\}, \{1, 3\}, \{2, 3\}\})*\{3\}|_{\{1, 2\}} = (\{1, 2\}, \{\O, \{1\}, \{2\}\}) \cong D_1$.

For the third, $(\{1, 2, 3\}, \{\O, \{2\}, \{3\}, \{1, 2\}, \{1, 3\}, \{1, 2, 3\}\})|_{\{2, 3\}} = (\{2, 3\}, \{\O, \{2\}, \{3\}\}) \cong D_1$.

For the fourth, $(\{1, 2, 3, 4\}, \{\O, \{1, 2\}, \{1, 3\}, \{1, 4\}, \{2, 3\}, \{2, 4\}, \{3, 4\}\})*\{1, 4\}|_{\{1, 2, 3\}} = (\{1, 2, 3\}, \{\O, \{1, 2\}, \{1, 3\}\}) \cong D_2$.

For the fifth, $(\{1,\, 2,\, 3,\, 4\},\; \{\O,\; \{1,\, 2\},\; \{1,\, 4\},\; \{2,\, 3\},\; \{3,\, 4\},\; \{1,\, 2,\, 3,\, 4\}\})|_{\{1, 2, 4\}} = (\{1, 2, 4\}, \{\O, \{1, 2\}, \{1, 4\}\}) \cong D_2$. \qedsymbol\\

The following lemma was proven by Yan and Jin within the proof of Theorem 12 of [18] for $D_1$ and within the proof of Theorem 20 of [17] for $D_2$ with only slightly stronger assumptions on $D$, but is reproven here for completeness.\\

\noindent
\textbf{Lemma 16.} Suppose that $D = (E, \mathcal{F})$ is a normal delta-matroid and $^{\partial}w_D$ is a monomial. Then no restriction of $D$ is isomorphic to $D_1$ or $D_2$.\\

\noindent
\textbf{Proof.} Let $m = r(D_{\max}) = w(D)$.

First suppose that some restriction of $D$ is isomorphic to $D_1$, so that for some distinct $e, f \in E$, $D|_{\{e, f\}} = (\{e, f\}, \{\O, \{e\}, \{f\}\})$. If some $F \in \mathcal{F}_{\max}$ does not contain $e$, then we would have $r((D*e)_{\max}) = m+1$ and $r((D*e)_{\min}) = 0$, so $w(D*e) = w(D)+1$. However, this contradicts the assumption that $^{\partial}w_D$ is a monomial, and so every $F \in \mathcal{F}_{\max}$ contains $e$, and similarly, contains $f$. For any $F \in \mathcal{F}_{\max-1}$, since $F$ is a subset of some set in $\mathcal{F}_{\max}$ by Proposition 13, this implies that $F$ must contain at least one of $e$ or $f$. Thus $r((D*\{e, f\})_{\max}) \leq m$, but since $r((D*\{e, f\})_{\min}) = 1$, this implies that $w(D*\{e, f\}) \leq w(D)-1$. This contradicts $^{\partial}w_D$ being a monomial, and so no restriction of $D$ is isomorphic to $D_1$.

Now suppose that some restriction of $D$ is isomorphic to $D_2$, so that for some distinct $e, f, g \in E$, $D|_{\{e, f, g\}} = (\{e, f, g\}, \{\O, \{e, f\}, \{e, g\}\})$. For any $h \neq f, g$ in $E$, $\{f, g, h\} \notin \mathcal{F}$, because else applying the symmetric exchange axiom to this set and the empty set would imply that at least one of $\{f\}$, $\{g\}$, or $\{f, g\}$ is in $\mathcal{F}$, but none of these are feasible. This shows that $r((D*\{f, g\})_{\min}) = 2$, and so since $w(D*\{f, g\}) = w(D)$, we must have $r((D*\{f, g\})_{\max}) = m+2$, meaning that for some $F \in \mathcal{F}_{\max}$, $F \cap \{f, g\} = \O$. This shows that $r_{D^*{}_{\min}}(\{e, f, g\}) \geq 2$. If there exists some $F \in \mathcal{F}_{\max}$ with $F \cap \{e, f\} = \O$, then $r((D*\{e, f\})_{\max}) = m+2$, but since $r((D*\{e, f\})_{\min}) = 0$, this would imply that $w(D*\{e, f\}) = w(D)+2$, a contradiction. Thus any $F \in \mathcal{F}_{\max}$ contains at least one of $e$ or $f$, and so $r_{D^*{}_{\min}}(\{e, f\}) \leq 1$. Similarly, $r_{D^*{}_{\min}}(\{e, g\}) \leq 1$. Since $r((D*e)_{\min}) = 1$ and $w(D*e) = w(D)$, we must have $r((D*e)_{\max}) = m+1$, meaning that some $F \in \mathcal{F}_{\max}$ does not contain $e$ and thus $r_{D^*{}_{\min}}(\{e\}) = 1$. Since the rank function of a matroid is submodular, $3 \leq r_{D^*{}_{\min}}(\{e, f, g\}) + r_{D^*{}_{\min}}(\{e\}) = r_{D^*{}_{\min}}(\{e, f\} \cup \{e, g\}) + r_{D^*{}_{\min}}(\{e, f\} \cap \{e, g\}) \leq r_{D^*{}_{\min}}(\{e, f\}) + r_{D^*{}_{\min}}(\{e, g\}) \leq 2$. This is clearly a contradiction, and so no restriction of $D$ is isomorphic to $D_2$. \qedsymbol\\

\noindent
\textbf{Proof of Theorem 2.} Suppose that the claim is false, and choose some nonbinary delta-matroid $D = (E, \mathcal{F})$ such that $^{\partial}w_D$ is a monomial. By Lemma 15, $D$ must have a minor isomorphic to $D_i$ for some $i \in \{1, 2\}$. Since $D_i$ is normal, Lemma 12 shows that there exist $F \in \mathcal{F}$ and $A \subseteq E$ such that $D_i \cong (D*F) \backslash A = (D*F)|_{A^c}$. Then since twisting by a set has no effect on the twist polynomial, $^{\partial}w_{D*F}$ is a monomial. However, $D*F$ is normal and has a restriction isomorphic to $D_i$, which contradicts Lemma 16. \qedsymbol\\

\noindent
\textbf{Acknowledgements}\\

The author would like to thank Leon Lozinskiy for helpful discussions and computational work, and Sergei Chmutov for valuable feedback on the paper and guidance throughout the project. The author is also thankful to the anonymous referees, whose useful comments helped improve the paper. This work has been done as a part of the undergraduate research program “Knots and Graphs” at the Ohio State University, during the summer of 2022. The author is grateful to the OSU Honors Program Research Fund for financial support.\\

\noindent
\textbf{References}\\

\begin{hangparas}{1.5em}{1}
[1] J. E. Bonin, C. Chun and S. D. Noble, Delta-matroids as subsystems of sequences of Higgs lifts, \textit{Adv. in Appl. Math}. \textbf{126} (2021) 101910.\\
\end{hangparas}

\begin{hangparas}{1.5em}{1}
[2] A. Bouchet, Greedy algorithm and symmetric matroids, \textit{Math. Program}. \textbf{38} (1987) 147–159.\\
\end{hangparas}

\begin{hangparas}{1.5em}{1}
[3] A. Bouchet, Representability of $\Delta$-matroids, in \textit{Colloq. Math. Soc. János Bolyai 52, Combinatorics}, Eger, Hungary, 1987, pp. 167-182.\\
\end{hangparas}

\begin{hangparas}{1.5em}{1}
[4] A. Bouchet, Maps and delta-matroids, \textit{Discrete Math}. \textbf{78} (1989) 59–71.\\
\end{hangparas}

\begin{hangparas}{1.5em}{1}
[5] A. Bouchet and A. Duchamp, Representability of $\Delta$-matroids over \textit{GF}(2), \textit{Linear Algebra Appl}. \textbf{146} (1991) 67–78.\\
\end{hangparas}

\begin{hangparas}{1.5em}{1}
[6] R. Brijder and H. J. Hoogeboom, Interlace polynomials for multimatroids and delta-matroids, \textit{European J. Combin}. \textbf{40} (2014) 142-167.\\
\end{hangparas}

\begin{hangparas}{1.5em}{1}
[7] S. Chmutov, Generalized duality for graphs on surfaces and the signed Bollobás-Riordan polynomial, \textit{J. Combin. Theory Ser. B}. \textbf{99} (2009) 617–638.\\
\end{hangparas}

\begin{hangparas}{1.5em}{1}
[8] S. Chmutov and F. Vignes-Tourneret, On a conjecture of Gross, Mansour and Tucker, \textit{European J. Combin}. \textbf{97} (2021) 103368.\\
\end{hangparas}

\begin{hangparas}{1.5em}{1}
[9] C. Chun, I. Moffatt, S. D. Noble and R. Rueckriemen, Matroids, delta-matroids and embedded graphs, \textit{J. Combin. Theory Ser. A}. \textbf{167} (2019) 7–59.\\
\end{hangparas}

\begin{hangparas}{1.5em}{1}
[10] D. Funk, D. Mayhew and S. D. Noble, How many delta-matroids are there?, \textit{European J. Combin}. \textbf{69} (2018) 149-158.\\
\end{hangparas}

\begin{hangparas}{1.5em}{1}
[11] J. F. Geelen, S. Iwata and K. Murota, The linear delta-matroid parity problem, \textit{J. Combin. Theory Ser. B}. \textbf{88} (2003) 377–398.\\
\end{hangparas}

\begin{hangparas}{1.5em}{1}
[12] J. L. Gross, T. Mansour and T. W. Tucker, Partial duality for ribbon graphs, I: Distributions, \textit{European J. Combin}. \textbf{86} (2020) 103084.\\
\end{hangparas}

\begin{hangparas}{1.5em}{1}
[13] I. Moffatt, Delta-matroids for graph theorists, in \textit{London Math. Soc. Lecture Note Ser. 456}, Cambridge Univ. Press, Cambridge, 2019, pp. 167–220.\\
\end{hangparas}

\begin{hangparas}{1.5em}{1}
[14] J. Oxley, Matroid theory, 2nd edn, Oxford University Press, New York, 2011.\\
\end{hangparas}

\begin{hangparas}{1.5em}{1}
[15] Q. Yan and X. Jin, Counterexamples to a conjecture by Gross, Mansour and Tucker on partial-dual genus polynomials of ribbon graphs, \textit{European J. Combin}. \textbf{93} (2021) 103285.\\
\end{hangparas}

\begin{hangparas}{1.5em}{1}
[16] Q. Yan and X. Jin, Partial-dual genus polynomials and signed intersection graphs, Preprint arXiv:2102.01823v1 [math.CO].\\
\end{hangparas}

\begin{hangparas}{1.5em}{1}
[17] Q. Yan and X. Jin, Twist polynomials of delta-matroids, \textit{Adv. in Appl. Math}. \textbf{139} (2022) 102363.\\
\end{hangparas}

\begin{hangparas}{1.5em}{1}
[18] Q. Yan and X. Jin, Twist monomials of binary delta-matroids, Preprint arXiv:2205.03487v1 [math.CO].\\
\end{hangparas}

\end{document}